\documentclass{article}

\begin{document}

\begin{center}
\textbf{On hypercomplexifying real forms of arbitrary rank}

Dennis Hou

\textit{dhou@uci.edu}

Abstract

\end{center}
For certain problems involving vector fields, it is possible 
to find an associated imaginary field that, in conjunction with 
the first, forms a complex field for which the equation can be 
solved. This result is generalized to arbitrary Clifford algebras, 
followed by quaternionic vectors as a special case. All results 
are shown to reduce to the established method of complexifying 
vector fields. For simplicity, differential forms are used rather 
than vector notation.


\section{Introduction}

The major feature of the algorithm in [1] for solving the Trkalian 
equation is writing a vector field as the real part of a complex 
field. Specifically, given the Monge potentials f, g, and h, 
a vector field v can be written as follows:

\begin{center}
v = f\ensuremath{\nabla}g + \ensuremath{\nabla}h.

\end{center}
Here, v is the real part of a vector field c = e$^{ig}$\ensuremath{\nabla}F, 
where F is defined by

\begin{center}
F \ensuremath{=} (h + if)e$^{-ig}$.

\end{center}
It follows that the imaginary part, w, the ``topological dual'' 
(henceforth ``Baldwin dual'') of v, is given by

\begin{center}
w = -h\ensuremath{\nabla}g + \ensuremath{\nabla}f.

\end{center}
In terms of differential forms [2],

\begin{center}
A = f dg + dh,
A$^{\ensuremath{\sim}}$ = -h dg + df.

\end{center}
It is possible to generalize this to higher dimensions. For instance, 
a vector field of five to eight dimensions can be thought of 
as the real part of a quaternionic field.

\section{Generalized Baldwin Dual}

Consider the one-form A of rank 2n - 1, where n = 2$^{k-1}$ for 
some integer k. The extended non-unique Clebsch decomposition is

\begin{center}
A = f$_{1}$df$_{2}$ + f$_{3}$df$_{4}$ +... + df$_{2n-1}$.

\end{center}
Then define F such that

\begin{center}
F \ensuremath{=} (f$_{2n-1}$ + e$_{1}$f$_{1}$ + e$_{2}$f$_{3}$ +... + e$_{n-1}$f$_{2n-3}$)exp(-e$_{1}$f$_{2}$ 
- e$_{2}$f$_{4}$ -... - e$_{n-1}$f$_{2n-2}$).

\end{center}

\begin{center}
\={F} \ensuremath{=} (f$_{2n-1}$ - e$_{1}$f$_{1}$ - e$_{2}$f$_{3}$ -... - e$_{n-1}$f$_{2n-3}$)exp(e$_{1}$f$_{2}$ 
+ e$_{2}$f$_{4}$ +... + e$_{n-1}$f$_{2n-2}$),

dF = ((df$_{2n-1}$ + e$_{1}$df$_{1}$ + e$_{2}$df$_{3}$ +... + e$_{n-1}$df$_{2n-3}$) - 
(f$_{2n-1}$ + e$_{1}$f$_{1}$ + e$_{2}$f$_{3}$ +... + e$_{n-1}$f$_{2n-3}$)(e$_{1}$df$_{2}$ + 
e$_{2}$df$_{4}$ +...\\
+ e$_{n-1}$f$_{2n-2}$))exp(-e$_{1}$f$_{2}$ - e$_{2}$f$_{4}$ -... - e$_{n-1}$f$_{2n-2}$),

d\={F} = ((df$_{2}$n$_{-1}$ - e$_{1}$df$_{1}$ - e$_{2}$df$_{3}$ -... - e$_{n-1}$df$_{2n-3}$) 
+ (f$_{2}$n$_{-1}$ - e$_{1}$f$_{1}$ - e$_{2}$f$_{3}$ -... - e$_{n-1}$f$_{2n-3}$)(e$_{1}$df$_{2}$ 
+ e$_{2}$df$_{4}$ \\
+... + e$_{n-1}$f$_{2n-2}$))exp(e$_{1}$f$_{2}$ + e$_{2}$f$_{4}$ +... + e$_{n-1}$f$_{2n-2}$).

\end{center}
This immediately leads to a generalization of the main result, 
namely that A can be represented as the real part of a one-form 
acting upon an n-dimensional Clifford algebra. Specifically,

\begin{center}
A = (exp(e$_{1}$f$_{2}$ + e$_{2}$f$_{4}$ +... + e$_{n-1}$f$_{2n-2}$) dF + exp(-e$_{1}$f$_{2}$ 
- e$_{2}$f$_{4}$ -... - e$_{n-1}$f$_{2n-2}$) d\={F})/2,

\end{center}
which is equivalent to Re(exp(e$_{1}$f$_{2}$ + e$_{2}$f$_{4}$ +... + e$_{n-1}$f$_{2n-2}$) 
dF). Defining the dual as the entire non-real portion fails because 
of the presence of multiple bases, so let A$^{\ensuremath{\sim}k}$ denote 
the Baldwin dual of A that corresponds to the kth base (excluding 
unity) on the algebra associated with the rank of the form. Extracting 
from exp(e$_{1}$f$_{2}$ + e$_{2}$f$_{4}$ +... + e$_{2}$n-1$_{-1}$f$_{2}$n$_{-2}$) dF,

\begin{center}
A$^{\ensuremath{\sim}1}$ = df$_{1}$ - f$_{2n-1}$df$_{2}$ - f$_{3}$df$_{6}$ + f$_{5}$df$_{4}$ - f$_{7}$df$_{10}$ 
+ f$_{9}$df$_{8}$ - \dots  - f$_{2n-5}$df$_{2n-2}$ + f$_{2n-3}$df$_{2n-4}$,

A$^{\ensuremath{\sim}2}$ = df$_{3}$ - f$_{2n-1}$df$_{4}$ + f$_{1}$df$_{6}$ - f$_{5}$df$_{2}$ + f$_{7}$df$_{12}$ 
- f$_{9}$df$_{14}$ +\dots  + f$_{2n-5}$df$_{2n-8}$ + - f$_{2n-3}$df$_{2n-6,}$

\end{center}
and so on. For n = 2, the imaginary dual reduces to the known 
form. There is also a natural extension of the set of conditions 
that corresponds to the known criteria for Baldwin duality over 
a complex field with respect to g,

\begin{center}
dA + A$^{\ensuremath{\sim}}$ dg = 0,

A$^{\ensuremath{\sim}}$ \ensuremath{\wedge} dA = 0,

A \ensuremath{\wedge} dA$^{\ensuremath{\sim}}$ = 0,

A \ensuremath{\wedge} dA = A$^{\ensuremath{\sim}}$ \ensuremath{\wedge} dA$^{\ensuremath{\sim}}$.

\end{center}
Specifically,

\begin{center}
dA + f$_{2}$A$^{\ensuremath{\sim}1}$ + f$_{4}$A$^{\ensuremath{\sim}2}$ + f$_{6}$A$^{\ensuremath{\sim}3}$ +... 
f$_{2n-2}$A$^{\ensuremath{\sim}n-1}$ = 0.

\end{center}
However, it is difficult to write down a generalization for

\begin{center}
dA + A$^{\ensuremath{\sim}}$ dg = 0.

\end{center}
Given the value of n, though, generating equations for the duals 
becomes simple. Using

\begin{center}
d(exp(e$_{1}$f$_{2}$ + e$_{2}$f$_{4}$ +... + e$_{n-1}$f$_{2n-2}$) dF) = dA + e$_{1}$dA\ensuremath{\sim}1 
+ e$_{2}$dA\ensuremath{\sim}2 +... + e$_{n-1}$ dA\ensuremath{\sim}n-1 = (e$_{1}$f$_{2}$ + e$_{2}$f$_{4}$ 
+... + e$_{n-}$f$_{2n-2}$)(exp(e$_{1}$f$_{2}$ + e$_{2}$f$_{4}$ +... + e$_{n-1}$f$_{2n-2}$) 
dF)

\end{center}
and equating parts, the analogues quickly arise. Restraints on 
the duals can be determined from this:

\begin{center}
f$_{2}$dA$^{\ensuremath{\sim}1}$ + f$_{4}$dA$^{\ensuremath{\sim}2}$ + f$_{6}$A$^{\ensuremath{\sim}3}$ +... 
f$_{2n}$dA$^{\ensuremath{\sim}n-1}$ = 0,

\end{center}

\begin{center}
dA$^{\ensuremath{\sim}k}$ \ensuremath{\wedge} A \ensuremath{\wedge} A$^{\ensuremath{\sim}1}$ \ensuremath{\wedge} 
\ensuremath{A}$^{\ensuremath{\sim}2}$... \ensuremath{\wedge} \ensuremath{A}$^{\ensuremath{\sim}k-1}$ \ensuremath{\wedge} 
\ensuremath{A}$^{\ensuremath{\sim}k+1}$ \ensuremath{\wedge}... \ensuremath{\wedge} \ensuremath{A}$^{\ensuremath{\sim}n 
-1}$ \ensuremath{=} 0,

A \ensuremath{\wedge} A$^{\ensuremath{\sim}1}$ \ensuremath{\wedge} \ensuremath{A}$^{\ensuremath{\sim}2}$... \ensuremath{\wedge} 
\ensuremath{A}$^{\ensuremath{\sim}n -1}$ \ensuremath{\neq} 0,

A \ensuremath{\wedge} dA$^{n-1}$ = A$^{\ensuremath{\sim}1}$ \ensuremath{\wedge} (dA$^{\ensuremath{\sim}1}$)$^{n-1}$ 
= A$^{\ensuremath{\sim}2}$ \ensuremath{\wedge} (dA$^{\ensuremath{\sim}2}$)$^{n-1}$ =... = A$^{\ensuremath{\sim}n-1}$ \ensuremath{\wedge} 
(dA$^{\ensuremath{\sim}n-1}$)$^{n-1}$

\end{center}
\section{Special Case: Quaternions}

The above generalization finds the Baldwin dual for any form 
of rank 2n - 1 to a n-dimensional algebra. Forms of rank 2n can 
be easily constructed simply by making the modulus of the hypercomplex 
form something other than unity. For other ranks, there at first 
seems to be no standard method of construction. Consider the 
special case n = 3:

\begin{center}
A = f$_{1}$df$_{2}$ + f$_{3}$df$_{4}$ + f$_{5}$df$_{6}$ + df$_{7}$,

F \ensuremath{=} (f$_{7}$ + if$_{1}$ + jf$_{3}$ + kf$_{5}$)exp(-if$_{2}$ - jf$_{4}$ - kf$_{6}$),

\={F} = (f$_{7}$ - if$_{1}$ - jf$_{3}$ - kf$_{5}$)exp(if$_{2}$ + jf$_{4}$ + kf$_{6}$),

dF = (df$_{7}$ + i df$_{1}$ + j df$_{3}$ + k df$_{5}$)exp(-if$_{2}$ - jf$_{4}$ - kf$_{6}$) 
- (f$_{7}$ + if$_{1}$ + jf$_{3}$ + kf$_{5}$)(i df$_{2}$ + j df$_{4}$ + k df$_{6}$)exp(-if$_{2}$ 
- jf$_{4}$ - kf$_{6}$),

d\={F} = (df$_{7}$ - i df$_{1}$ - j df$_{3}$ - k df$_{5}$)exp(if$_{2}$ + jf$_{4}$ 
+ kf$_{6}$) + (f$_{7}$ - if$_{1}$ - jf$_{3}$ - kf$_{5}$)(i df$_{2}$ + j df$_{4}$ + k 
df$_{6}$)exp(if$_{2}$ + jf$_{4}$ + kf$_{6}$),

A = Re(exp(if$_{2}$ + jf$_{4}$ + kf$_{6}$) dF).

\end{center}
The Baldwin duals are as follows:

\begin{center}
A$^{\ensuremath{\sim}1}$ = df$_{1}$ - f$_{7}$df$_{2}$ - f$_{3}$df$_{6}$ + f$_{5}$df$_{4}$,

\end{center}

\begin{center}
A$^{\ensuremath{\sim}2}$ = df$_{3}$ - f$_{7}$df$_{4}$ + f$_{1}$df$_{6}$ - f$_{5}$df$_{2}$,

A$^{\ensuremath{\sim}3}$ = df$_{5}$ - f$_{7}$df$_{6}$ - f$_{1}$df$_{4}$ + f$_{3}$df$_{2}$.

\end{center}
Now let B = A + iA$^{\ensuremath{\sim}1}$ + jA$^{\ensuremath{\sim}2}$ + kA$^{\ensuremath{\sim}3}$. 
Then

\begin{center}
dB = (if$_{2}$ + jf$_{4}$ + kf$_{6}$) B,

\end{center}
which leads to,

\begin{center}
dA + f$_{2}$A$^{\ensuremath{\sim}1}$ + f$_{4}$A$^{\ensuremath{\sim}2}$ + f$_{6}$A$^{\ensuremath{\sim}3}$ = 
0,

dA$^{\ensuremath{\sim}1}$ - f$_{2}$A + f$_{6}$A$^{\ensuremath{\sim}2}$ - f$_{4}$A$^{\ensuremath{\sim}3}$ = 
0,

dA$^{\ensuremath{\sim}2}$ - f$_{4}$A - f$_{6}$A$^{\ensuremath{\sim}1}$ + f$_{2}$A$^{\ensuremath{\sim}3}$ = 
0,

dA$^{\ensuremath{\sim}3}$ - f$_{6}$A + f$_{4}$A$^{\ensuremath{\sim}1}$ - f$_{2}$A$^{\ensuremath{\sim}2}$ = 
0.

\end{center}
Using the exterior derivative and product on this set of equations 
leads to more conditions satisfied by A and its Baldwin duals.

\begin{center}
f$_{2}$dA$^{\ensuremath{\sim}1}$ + f$_{4}$dA$^{\ensuremath{\sim}2}$ + f$_{6}$dA$^{\ensuremath{\sim}3}$ = 0,

f$_{2}$dA - f$_{6}$dA$^{\ensuremath{\sim}2}$ + f$_{4}$dA$^{\ensuremath{\sim}3}$ = 0,

f$_{4}$dA + f$_{6}$dA$^{\ensuremath{\sim}1}$ - f$_{2}$dA$^{\ensuremath{\sim}3}$ = 0,

f$_{6}$dA - f$_{4}$dA$^{\ensuremath{\sim}1}$ + f$_{2}$dA$^{\ensuremath{\sim}2}$ = 0,

dA \ensuremath{\wedge} A$^{\ensuremath{\sim}1}$ \ensuremath{\wedge} A$^{\ensuremath{\sim}2}$ \ensuremath{\wedge} A$^{\ensuremath{\sim}3}$ 
= 0,

dA$^{\ensuremath{\sim}1}$ \ensuremath{\wedge} A \ensuremath{\wedge} A$^{\ensuremath{\sim}2}$ \ensuremath{\wedge} A$^{\ensuremath{\sim}3}$ 
= 0,

dA$^{\ensuremath{\sim}2}$ \ensuremath{\wedge} A \ensuremath{\wedge} A$^{\ensuremath{\sim}1}$ \ensuremath{\wedge} 
A$^{\ensuremath{\sim}3}$ = 0,

dA$^{\ensuremath{\sim}3}$ \ensuremath{\wedge} A \ensuremath{\wedge} A$^{\ensuremath{\sim}1}$ \ensuremath{\wedge} A$^{\ensuremath{\sim}2}$ 
= 0,

\end{center}

\begin{center}
A \ensuremath{\wedge} dA$^{3}$ = A$^{\ensuremath{\sim}1}$ \ensuremath{\wedge} (dA$^{\ensuremath{\sim}1}$)$^{3}$ 
= A$^{\ensuremath{\sim}2}$ \ensuremath{\wedge} (dA$^{\ensuremath{\sim}2}$)$^{3}$ = A$^{\ensuremath{\sim}3}$ \ensuremath{\wedge} 
(dA$^{\ensuremath{\sim}3}$)$^{3}$.

\end{center}
Special cases that reduce the rank of A to six are not interesting, 
because they can only arise by eliminating one of the zero-form 
potentials, (or else the rank would drop to five,) but none of 
the constraints rely on these. However, for f$_{5}$ = f$_{6}$ = 0, the 
equations are reduced as follows:

\begin{center}
dA + f$_{2}$A$^{\ensuremath{\sim}1}$ + f$_{4}$A$^{\ensuremath{\sim}2}$ = 0,

dA$^{\ensuremath{\sim}1}$ - f$_{2}$A - f$_{4}$A$^{\ensuremath{\sim}3}$ = 0,

dA$^{\ensuremath{\sim}2}$ - f$_{4}$A + f$_{2}$A$^{\ensuremath{\sim}3}$ = 0,

dA$^{\ensuremath{\sim}3}$ + f$_{4}$A$^{\ensuremath{\sim}1}$ - f$_{2}$A$^{\ensuremath{\sim}2}$ = 0.

\end{center}

\begin{center}
f$_{2}$dA$^{\ensuremath{\sim}1}$ + f$_{4}$dA$^{\ensuremath{\sim}2}$ = 0,

f$_{2}$dA + f$_{4}$dA$^{\ensuremath{\sim}3}$ = 0,

f$_{4}$dA - f$_{2}$dA$^{\ensuremath{\sim}3}$ = 0,

f$_{4}$dA$^{\ensuremath{\sim}1}$ - f$_{2}$dA$^{\ensuremath{\sim}2}$ = 0,

dA \ensuremath{\wedge} A$^{\ensuremath{\sim}1}$ \ensuremath{\wedge} A$^{\ensuremath{\sim}2}$ = 0,

dA$^{\ensuremath{\sim}1}$ \ensuremath{\wedge} A \ensuremath{\wedge} A$^{\ensuremath{\sim}3}$ = 0,

dA$^{\ensuremath{\sim}2}$ \ensuremath{\wedge} A \ensuremath{\wedge} A$^{\ensuremath{\sim}3}$ = 0,

dA$^{\ensuremath{\sim}3}$ \ensuremath{\wedge} A$^{\ensuremath{\sim}1}$ \ensuremath{\wedge} A$^{\ensuremath{\sim}2}$ 
= 0,

A \ensuremath{\wedge} dA$^{2}$ = A$^{\ensuremath{\sim}1}$ \ensuremath{\wedge} (dA$^{\ensuremath{\sim}1}$)$^{2}$ 
= A$^{\ensuremath{\sim}2}$ \ensuremath{\wedge} (dA$^{\ensuremath{\sim}2}$)$^{2}$ = A$^{\ensuremath{\sim}3}$ \ensuremath{\wedge} 
(dA$^{\ensuremath{\sim}3}$)$^{2}$.

\end{center}
For rank-three forms,

\begin{center}
f$_{3}$ = f$_{4}$ = f$_{5}$ = f$_{6}$ = 0,

\end{center}

\begin{center}
A$^{\ensuremath{\sim}2}$ = A$^{\ensuremath{\sim}3}$ = 0.

\end{center}
Furthermore, most of the constraints vanish, leaving

\begin{center}
dA + f$_{2}$A$^{\ensuremath{\sim}1}$ = 0,

dA$^{\ensuremath{\sim}1}$ - f$_{2}$A = 0,

dA \ensuremath{\wedge} A$^{\ensuremath{\sim}1}$ = 0,

dA$^{\ensuremath{\sim}1}$ \ensuremath{\wedge} A = 0,

A \ensuremath{\wedge} dA = A$^{\ensuremath{\sim}1}$ \ensuremath{\wedge} dA$^{\ensuremath{\sim}1}$,

\end{center}
which are precisely the equations specified in [1].\\
\section{Acknowledgments}

This paper is dedicated to Charlotte Hwa, of whom it is written\footnote{Proverbs 
xxvii, 17.}, ``Iron sharpeneth iron; so a man sharpeneth the countenance 
of his friend.''

\section{References}

\lbrack1\rbrack P. R. Baldwin and G. M. Townsend, \textit{Complex Trkalian Fields and
Solutions to Euler's Equations for the Ideal Fluid}, \textbf{Phys. Rev. E} 51
(1995) 2059-2068.\\
\lbrack2\rbrack H. Flanders, {\underline{Differential Forms with Applications
to the Physical Sciences}}, (Dover, N. Y., 1963).\\
\lbrack3\rbrack P. Lounesto, {\underline{Clifford Algebras and Spinors}},
(Cambridge University Press, Cambridge, 1997).\\

\end{document}